\pdfoutput=1
\documentclass[10pt,twocolumn,twoside]{IEEEtran}

\usepackage{packages}

\usepackage{commands}

\newcommand{\dict}[0]{\m{\Gamma}} 
\newcommand{\dictv}[0]{\vec{\gamma}} 
\newcommand{\netstate}[0]{\Xi} 
\newcommand{\stateset}[0]{\mathbb{B}} 
\newcommand{\edit}[1]{#1}

\newcommand{\editbegin}{}
\newcommand{\editend}{}

\begin{document}


\begin{acronym}[DOTDMA]

\acro{API}{application programming interface} 
\acro{AWSS}{asymptotically wide sense stationary}
\acro{AWGN}{additive white gaussian noise}
\acro{BB}{block-based}
\acro{BIBO}{bounded input bounded output}
\acro{BP}{basis pursuit}
\acro{CLPPC}{closed loop packet predictive controller}
\acro{CDF}{cohen-daubechies-feauveau}
\acro{CPU}{central processing unit}
\acro{CUDA}{compute unified device architecture}
\acro{CU}{compute unit}
\acro{CWT}{continuous wavelet transform}
\acro{DWT}{discrete wavelet transform}
\acro{DPR}{dynamic partial reconfiguration}
\acro{DtoH}{device to host}
\acro{ECC}{error correcting code}
\acro{ECDQ}{entropy-coded dithered (lattice) quantizer}
\acro{FB}{functional block}
\acro{FLOP}{floating point operation}
\acro{FLOPS}{floating point operations per second}
\acro{FPGA}{field programmable gate array}
\acro{GPU}{graphics processing unit}
\acro{GR}{Gaussian with refinements}
\acro{GSR}{Gaussian with scaled refinements}
\acro{GR2}{Gaussian dictionaries with refinements for the two state system}
\acro{GSR2}{Gaussian dictionaries with scaled refinements for the two state system}
\acro{HtoD}{host to device}
\acro{HD}{high definition}
\acro{IID}[i.i.d.]{independent and identically distributed}
\acro{IP}{integer programming}
\acro{JPEG}{joint photographic experts group}
\acro{LTI}{linear time invariant}
\acro{LQR}{linear quadratic regulator}
\acro{LB}{line-based}
\acro{MPC}{model predictive control}
\acro{MSE}{mean squared error}
\acro{MIMD}{multiple instruction multiple data}
\acro{MP}{matching pursuit}
\acro{MJLS}{Markov jump linear system}
\acro{MSS}{mean square stable}
\acro{NCS}{networked control system}
\acro{OMP}{orthogonal matching pursuit}
\acro{OpenCL}{open computing language}
\acro{PDF}{probability density function}
\acro{PR}{partial reconfigurable}
\acro{PRR}{partial reconfigurable region}
\acro{PRM}{partial reconfigurable module}
\acro{PE}{processing element}
\acro{PSNR}{peak signal to noise ratio}
\acro{PSD}{power spectral density}
\acro{PPC}{packetized predictive control}
\acro{RC}{row-column}
\acro{REC}[RC]{Reconfigurable}
\acro{SNR}{signal to noise ratio}
\acro{SIMD}{single instruction multiple data}
\acro{SPARC}{sparse regression codes}
\acro{VQ}{vector quantizer}

\end{acronym}

%
\title{Shaped Gaussian Dictionaries for Quantized Networked Control Systems with Correlated Dropouts}
%
%
%

\author{Edwin G.W.~Peters, Daniel E.~Quevedo,~\IEEEmembership{Senior Member,~IEEE,} Jan Østergaard,~\IEEEmembership{Senior Member,~IEEE}
\thanks{E. G. W. Peters is with the School of Electrical Engineering \& Computer Science, The University of Newcastle, Callaghan NSW 2308 Australia, {\tt edwin.g.w.peters@gmail.com}.}
\thanks{D. E. Quevedo is with the Department of Electrical Engineering (EIM-E) at the University of Paderborn, 33098 Paderborn, Germany {\tt dquevedo@ieee.org}.}
\thanks{J. Østergaard is with the Department of Electronic Systems, Aalborg University, Aalborg 9220, Denmark {\tt janoe@ieee.org}.}
\thanks{The work of J. {\O}stergaard is financially supported by VILLUM FONDEN Young Investigator Programme, Project No. 10095.}
}

\markboth{}%
{Peters \MakeLowercase{\textit{et al.}}: Shaped Gaussian Dictionaries for Quantized Networked Control Systems with Correlated Dropouts }

\maketitle

\begin{abstract}
This paper studies fixed rate vector quantisation for noisy \acp{NCS} with correlated packet dropouts. In particular, a discrete-time \acl{LTI} system is to be controlled over an error-prone digital channel. The controller uses (quantized) \acl{PPC} to reduce the impact of packet losses. The proposed \acl{VQ} is based on \ac{SPARC}, which have recently been shown to be efficient in open-loop systems when coding white Gaussian sources. The dictionaries in existing design of \acp{SPARC} consist of \ac{IID} Gaussian entries. However, we show that a significant gain can be achieved by using Gaussian dictionaries that are shaped according to the second-order statistics of the \ac{NCS} in question. Furthermore, to avoid training of the dictionaries, we provide closed-form expressions for the required second-order statistics in the absence of quantization. 
\end{abstract}


\acresetall


\section{Introduction}
\label{cha:intr}

\IEEEPARstart{L}{inear} Time Invariant (\acs{LTI}) control systems are today used in many different places and
situations. These systems all have in common, that there is feedback from the
system (plant) to be controlled to the controller, in order to maintain the plant in a
desired state. In some cases it might be desired to have the controller at
a different physical location than the plant, in which case a wired or wireless network connects
them. This topology is called a \ac{NCS} and can have many advantages such as lower cost, higher reliability and easier maintenance. However, other challenges arise, including bit rate limitations, random delays and breakdowns, which leave the plant in open loop operation and can have severe consequences depending on the situation \cite{4118465,4626109,5345802,5484486,4814536}.

\acp{NCS} with random \ac{IID} packet dropouts and delays have been introduced and analyzed in previous works under the terms ``\acl{PPC}'' (\acsu{PPC}) and ``receding horizon control'', where plant input predictions are transmitted and stored in a buffer. These can then be used in case a packet that contains input signals to the plant does not arrive. 
In particular, the paper \cite{5648448} shows stability results for cases where the maximum number of consecutive packet dropouts is bounded, whereas \cite{5471054} investigates \acp{NCS} with bounded time-delays. The works \cite{Quevedo20121803,5742978} analyze mean square- and stochastic stability of \acp{NCS} based on a Markov dropout model where an unbounded number of consecutive packet dropouts may occur.

Quantization within the \ac{NCS} has been studied e.g. in \cite{5742978,6137920,RNC:RNC887}, where \cite{RNC:RNC887} forces the controller to select the control vector from a finite constrained set of vectors using a nearest neighbor \ac{VQ}, and analyses the closed loop behavior of these. In \cite{6137920}, an \ac{ECDQ} is used, and closed loop stability is investigated using linear matrix inequalities based on \acp{MJLS}. Optimal rates for the entropy coder are calculated based on the statistics of the \ac{NCS}. 
\ac{MJLS} stability using \ac{ECDQ} is investigated in \cite{5742978}, where the authors provide bounds on the maximum packet dropout rates which preserve stability in the system.
The work \cite{6426647} relates the \ac{PPC} to problems solved in compressed sensing, and investigates sparse representations of the control vector using different techniques including \acl{OMP}. The authors furthermore provide sufficient conditions for stability when the controller is used on a network where bounded packet dropouts occur.

The contributions to \ac{NCS} in this paper focus on using a fixed rate \ac{VQ} in \ac{PPC}. \edit{\acp{VQ}, in comparison to scalar quantizers, have the ability to operate on multiple dimensions. This opens up for concepts that are not applicable on scalars, and can at least match the performance of a scalar quantizer, although it often does better \cite[Theorem 10.1.1]{Gersho:1991:VQS:128857}. Fixed rate \acp{VQ} further have the advantage, that the network bandwidth requirements are fixed. Since most network protocols (such as the IEEE 802.15.4 based WirelessHART) utilize slotted transmissions, using a fixed rate \ac{VQ} one can fit the information in a slot. This has a clear advantage over variable rate coding approaches such as entropy based coding, where smaller packets can be obtained, but it is still necessary to reserve additional bandwidth in case a larger packet has to be transmitted. In this work we utilize a \ac{VQ} which allows one to compress an entire vector at once.} The dictionary \edit{used by the \ac{VQ}} is inspired by \ac{SPARC}, which is introduced in \cite{6284210} as efficient codes to compress memoryless Gaussian sources.
In \cite{Quevedo20121803,4982641} the packet dropouts in the network are considered to be correlated as described by \cite{Elliott1963,pr:81}. This produces a more real-world description of the network where disturbances on the network can occur for small periods of time where more frequent packet dropouts take place. In this work we will adopt this class of network models and exploit this in the design of the dictionaries for the \ac{NCS}.

Additionally, we will implement a fixed rate \ac{VQ} in \acp{NCS} using receding horizon controllers, where we investigate and apply different methods to construct the dictionaries without the need of training data. We furthermore investigate network models for correlated packet dropouts and design quantizers for these. This results in a more realistic representation of the network between the controller and buffer. We finally provide simulation studies of the quantizer in network setups featuring correlated and \ac{IID} dropouts.

\textit{Notation:} 
$\dictv_{i}$ denotes the $i^{\mathrm{th}}$ column vector in matrix $\dict$. \edit{Let a matrix $\m{A} \in \mathbb{R}^{m\times m}$ denote an $m\times m$ matrix and $\m[T][]{A}$ its transpose, while $\vec{B}\in \R^{m}$ is a vector containing $m$ elements. We further denote by $a_{j,i}$ the $j^{\mathrm{th}}$ element in the $i^{\mathrm{th}}$ column of a matrix $\m{A}$.} For any vector $\vec{x}$ and square matrix $\m{Q}$, we define $\|\vec{x}\|_{\bm{Q}}^{2} = \vec[T][]{x}\m{Q}\vec{x}$ and $\|\vec{x}\|^{2}=\vec[T][]{x}\vec{x}$. The $m\times m$ identity matrix is denoted by $\m[][m]{I}$ \edit{ and the $n\times m$ matrix containing zeros is denoted $\m{0}_{n\times m}$.}


\section{\ac{NCS} with fixed-rate quantizers}
\label{sec:acppc-with-fixed}

\begin{figure}[tp]
  \centering
  \def\svgwidth{0.9 \columnwidth}  
  \scriptsize{
  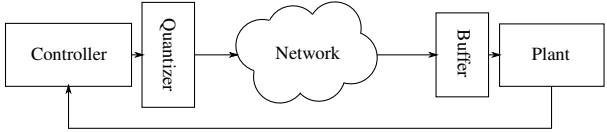
}
  \caption{An illustration of the \ac{NCS} considered with quantizer and buffer added. \edit{The controller computes a control signal that is quantized and afterwards transmitted over the network to the buffer. At every time step the buffer updates the actuator input to the plant.}}
  \label{fig:syssetupintro}
\end{figure}

The \ac{NCS} considered is shown in Figure \ref{fig:syssetupintro}, where the controller and the plant input are connected through a network in which packet loss occurs. In most previous work related to \acp{NCS}, the dropouts were assumed to be \ac{IID} \cite{6426647,5742978,6137920,5648448}, such that a packet dropout occurs with probability $p_{d}$ at every time-instance $k$. We utilize state feedback, where the feedback path throughout this work is assumed to be ideal with no packet dropouts.

\subsection{Packetized predictive control (\acs{PPC})}
\label{sec:acppc}
\editbegin
In this work we consider the system shown in \cref{fig:syssetupintro}, where the state of the plant is given by the recursion
\begin{align}\label{eq:nextstate}
  \vec[][]{x}(k+1) = \m{A}\vec[][]{x}(k) + \vec[][1]{B}u(k) + \vec[][2]{B}\omega(k), ~ k \in \mathbb{N}_{0},
\end{align}
with $\vec[][]{x}(k) \in \mathbb{R}^{p}$ describing the plant state at time $k$. In \cref{eq:nextstate} the system matrix of the plant is $\m{A} \in \mathbb{R}^{p\times p}$, $p$ is the dimension of the state vector, $u(k) \in \mathbb{R}$ is the control signal or input and $\vec[][1]{B} \in \mathbb{R}^{p}$. The plant is affected by a disturbance with zero-mean white Gaussian noise $\omega(k)$, that is applied through $\vec[][2]{B} \in \mathbb{R}^{p}$.

We only consider controllable systems in this work. This means, that the matrix
\begin{align}
  \mathcal{C} \triangleq \left[B,\,AB,\,A^{2}B,\,\dots,\,A^{p-1}B\right] \in \R^{p\times p}
\end{align}
has full row rank. This condition can be relaxed to saying that the system is stabilizable if the uncontrollable subspace of the system matrix $\m{A}$ has all of its eigenvalues strictly in the unit circle \cite{goodwin:2001:control}.

In this work, we assume that the probability for packet dropouts $p_{d}$ is non-zero which  means, that the computed control signal not always will be applied to the actuator. This has to be taken into account to achieve the desired closed loop performance. See for example \cite{RNC:RNC887,rawlings2009model} and other literature on \acf{PPC}.
The main idea in \ac{PPC} is that, based on current measurement data, the controller predicts which control signals will be applied to the actuators in the future. These control signals are transmitted over the network as a vector 
\begin{align}
\vec[][]{u}(k) = \left[u(k),\, u(k+1),\,\dots\,, u(k+N-1)\right]^{T} \in \mathbb{R}^{N},
\end{align}
where $N$ is the prediction horizon. The packets are received by a buffer, that stores the control signals. At every time instance, the buffer verifies whether the transmission of control data is successful. If a new packet with control signals arrives, the previous signal is replaced in the buffer and the first entry in $\vec[][]{u}(k)$ is applied to the actuator. In case a packet is lost, the buffer will apply the next control signal to the actuators. This leaves the operation of the buffer as follows
\begin{align}\label{eq:buffer}
  \vec{b}(k) = (1-d(k)) \vec[][]{u}(k) + d(k) \m{M} \vec[][]{b}(k-1) \in \R^{N},
\end{align}
where
\begin{align}
  \m{M} = \left[
    \begin{matrix}
      0 & 1 & 0 &\dots & 0\\
      \vdots & \ddots& \ddots & \ddots & \vdots\\
      0 & \dots & 0 & 1 & 0\\
      0 & \dots & 0 & 0 & 1\\      
      0 & \dots & 0 & 0 & 0
    \end{matrix}
  \right]
\end{align}
is the $N\times N$ matrix that cycles the buffer and the parameter $d(k)$ indicates whether a packet dropout occurred at time instance $k$. Here $d(k)=1$ when a packet dropout occurs. This happens with probability $p_{d}$. The actuator input $u(k)$ in \cref{eq:nextstate} is then taken from the top-most entry in the buffer, namely
\begin{align}
  u(k) = \vec{e}_{1} \vec{b}(k),
\end{align}
where the $1\times N$ vector $\vec{e}_{1}\triangleq [1,0,\dots,0]$.

In the \ac{PPC} formulation we use the linear quadratic cost function 
\begin{align}\label{eq:jleftvu-vkxright-=}
  \begin{split}
    &J\left(\vec{u}',\vec[][]{x}(k)\right) = \|\vec{x}'(N)\|^{2}_{{\m{X}}}
+\sum_{l=0}^{N-1}\left(\|\vec{x}'(l)\|_{{\m{Q}}}^{2}+Ru'(l)^{2}\right),
  \end{split}
\end{align}
where $\vec{u}'= \left[u'(0), \, u'(1),\,\dots,\,u'(N-1)\right]^{T}$ are the predicted actuator inputs and $\vec{x}'(k+l)$ are the predicted plant states for the inputs $u'(l)$ and are given by
\begin{align}
  \vec{x}'(l+1) = \m{A}\vec{x}'(l) + \vec{B}_{1}u'(l), \quad l \in \left\{0,1,\dots,N-1\right\}
\end{align}
with $\vec{x}'(0)= \vec{x}(k)$. 
Here $N\geq 1$ is the length of the prediction horizon which equals the size of the buffer in \cref{fig:syssetupintro}. The variables $\m{Q}\succeq 0 \in \R^{p\times p}$, $\m{X} \succeq 0 \in \R^{p\times p}$ and $R > 0 \in \R$ are weighting matrices and scalars that allow for trade-off between control performance and control effort \cite{goodwin:2001:control}. These parameters can be tuned in the design phase until the desired closed loop performance is obtained.

By defining the matrices
\begin{align}
  \m{\Phi} \triangleq \left[
  \begin{matrix}
    \vec[][1]{B} & \m{0}_{p} & \cdots & \m{0}_{p}\\
    \m{A}\vec[][1]{B} & \vec[][1]{B} & \cdots & \m{0}_{p}\\
    \vdots & \vdots & \ddots & \vdots \\
    \m{A}^{N-1}\vec[][1]{B} & \m{A}^{N-2}\vec[][1]{B} & \cdots & \vec[][1]{B}
  \end{matrix}
                                                                 \right] \in \R^{Np\times N}
\end{align}
and
\begin{align}
  \m{\Upsilon}  \triangleq \left[
  \begin{matrix}
    \m{A}\\
    \m{A}^{2}\\
    \vdots\\
    \m{A}^{N}
  \end{matrix}
  \right] \in \R^{Np\times p},
\end{align}
we can restate the cost function \cref{eq:jleftvu-vkxright-=} as 
\begin{align}\label{eq:closedloopcost2ppc}
  \begin{split}
    &J(\vec{u}',\vec{x}(k)) =  \vec[T]{x}(k) \m[T]{\Upsilon} \m{\Upsilon} \vec{x}(k) \mathbreak
 + \vec{u}'^{T} \m{W} \vec{u}'  + 2 \vec[T]{x}(k) \m{F} \vec{u}',
  \end{split}
\end{align}
where
\begin{align}
  &\m{W} \triangleq \bar{\m{R}} + \m{\Phi}^{T}\bar{\m{Q}}\m{\Phi}  \in \R^{N\times N}\\
  &\m{F} \triangleq \m{\Upsilon}^{T}\bar{\m{Q}}\m{\Phi}  \in \R^{p\times N}
\end{align}
and the weighting matrices are given by 
\begin{align}
  \begin{split}\label{eq:section2new:weights}
    &\m{\bar{Q}} \triangleq \operatorname{blockdiag}(\m{Q},\cdots,\m{Q},\m{X}) \in \R^{Np\times Np}\\
    &\m{\bar{R}} \triangleq \diag \left\{R, \cdots, R\right\} \in \R^{N\times N},
  \end{split}
\end{align}
where $\m{X} \in \R^{p\times p}$ is the symmetric positive definite solution to the discrete algebraic Riccati equation
\begin{align}\label{eq:dare}
  \begin{split}
    \m{X} =& \m[T]{A}\m{X}\m{A} + \m{Q} - \m[T]{A}\m{X}\vec[][1]{B}\left(R + \vec[T][1]{B}\m{X}\vec[][1]{B}\right)^{-1} \\
    &\vec[T][1]{B}\m{X}\m{A},
  \end{split}
\end{align}
which exists only if the system \cref{eq:nextstate} is stabilizable \cite{goodwin:2001:control}. The optimal control inputs are then found by
\begin{align}
  \vec{u}(k) \triangleq \argmin_{\vec{u}\in \R^{N}} J(\vec{u},\vec{x}(k)),
\end{align}
which can be solved analytically and results in the feedback law
\begin{align}
    \vec[][]{u}(k) = -\m{K}\vec[][]{x}(k),
\end{align}
where
\begin{align}\label{eq:section2new:ctrlGain}
  \m{K}= \m{W}^{\, -1} \m{F}^{T}
\end{align}
is the stabilizing feedback gain.

Since we in this work focus on bandwidth limited networks, it is of high interest that the size in bits of the control packets is small. This can be done using vector quantization. The papers \cite{6426647,5742978,6137920} investigated this idea using different quantizers and methods where variable-rate \acp{VQ} are used. Variable rate \acp{VQ} give the disadvantage, that the demands on the communication channel vary depending on the amount of bits that are required to store the control signal. This results in wasted resources on the channel in case the bit rate has high variations. 
To overcome this problem we, in the present work, propose to apply a fixed-rate \ac{VQ} in the \ac{NCS}, such that the requirements for the communication network do not vary. 
\editend

\subsection{Fixed rate quantizers for \ac{PPC}}
\label{sec:fixed-rate-quant}

\edit{In this work we utilize a \ac{VQ} to compress the control signal such that less network bandwidth is required to transmit the signal. A vector quantizer is a function, that maps a $N$-dimensional vector in a space $\R^{N}$ into a finite set $\mathbb{W}$ containing so-called codewords. This set is often referred to as the codebook, thus a vector quantizer $\mathcal{Q}: \R^{n} \rightarrow \mathbb{W}$ \cite{Gersho:1991:VQS:128857}.} 

We design the quantizer using a fixed dictionary $\dict$ that is split into $M$ sections which each contains $L$ codewords.  \edit{In this case, the codewords in the set $\mathbb{W}$ are composed of linear combinations of the codewords in the sections of the dictionary $\dict$.} The dictionary is known for both the quantizer and the buffer on the receiver side of the network. The dictionary design is heavily inspired by \ac{SPARC}, presented in \cite{6284210}. The structure of the dictionary in \ac{SPARC} can be used in the design of computationally efficient encoders \cite{DBLP:journals/corr/abs-1212-1707}. An additional feature of \ac{SPARC} is, that it has a low memory requirement, since only one section of the dictionary needs to be stored in the memory for every iteration. This can be advantageous when the algorithm is implemented on micro-controllers etc.
A key difference is, that while \ac{SPARC} is designed to compress large vectors containing \ac{IID} data, the designed dictionary will be used on shorter vectors containing memory due to the feedback in the control loop. Further details on the dictionary design are included in Section \ref{sec:dictionary-design}.
\begin{figure}[tp]
  \centering
  \def\svgwidth{0.8 \columnwidth}  
  \scriptsize{
 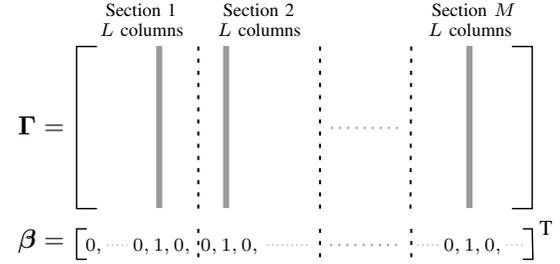
}
\caption{The proposed dictionary based on \ac{SPARC}, \edit{where the vector $\bm{\beta}$ only has one nonzero entry in each section $m$, that selects a vector from the dictionary $\dict$ and thereby can reconstruct the signal $u$ through the linear combination $\bar{\bm{u}} = \dict \bm{\beta}$.}
}
  \label{fig:sparc-dictionary}
\end{figure}

Figure \ref{fig:sparc-dictionary} illustrates the dictionary $\dict \in R^{N\times ML}$,
where $N$ is the horizon length of the controller. The vector $\m{\beta}$ is an $ML\times 1$ vector, which only contains one non-zero entry in every section $m \in \left\{1,2,\dots, M\right\}$, and this entry is $1$. A vector $\vec{\beta}$ can be used to estimate a signal $\vec{u}$ by
\begin{align}
  \vec{\bar{u}} = \dict \vec{\beta},
\end{align}
\edit{where 
  \begin{align}\label{eq:section2new:quant-error}
    \bar{\vec{u}} = \vec{u} - \epsilon,
  \end{align}
with $\epsilon$ being the quantization error.
}
Using this quantizer, only the vector $\vec{\beta}$ has to be transmitted to the receiver. This vector consists of $M \log_{2}(L)$ bits, which gives an effective bit rate of 
\begin{align}
  R = \frac{M \log_{2}(L)}{N} ~ \text{\SI{}{\bit\per\symbol}.}
\end{align}

From the set $\mathcal{B}$ containing all code words, we approximate $\vec{u}\in \mathbb{R}^{N}$ by finding the code word $\vec{\beta} \in \mathcal{B}$, for which $\dict \vec{\beta}$ is closest to $\vec{u}$. Thus, given a vector $\vec{u}$ and some fixed code book $\dict$, we have to solve the minimization problem
\editbegin
\begin{align}\label{eq:regressioncost}
  \vec[*]{\beta} = \argmin_{\vec{\beta}\in \mathcal{B}} \| \vec{u}- \dict\vec{\beta}\|^{2}.
\end{align}
The set $\mathcal{B}$ is non-convex, which means that the optimization (\ref{eq:regressioncost}) is a non-convex optimization problem. In fact, it is an NP-hard combinatorial problem \cite{6284210}. The optimization (\ref{eq:regressioncost}) needs to be solved on-line, which can be done using e.g. greedy algorithms such as \ac{MP} \cite{258082}. These algorithms do not necessarily provide the global optimal solution, but often result in a local minimum, which is sufficient depending on the application. The algorithm used in the present work is explained in Section \ref{sec:acmjls-networks-with}.

\begin{figure}[tp]
  \centering
  \def\svgwidth{\columnwidth}  
  \begin{scriptsize}
    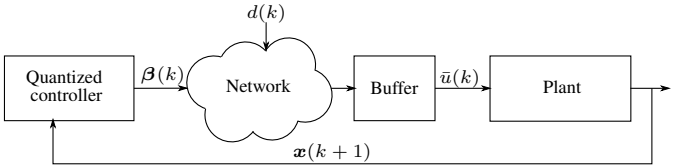
  \end{scriptsize}
  \caption{Closed loop \ac{PPC} with the quantizer integrated in the cost function according to (\ref{eq:closedloopcost2}).}
  \label{fig:clppc}
\end{figure}

We implement the quantizer into the controller as shown in Figure \ref{fig:clppc}, such that the controller directly operates on the finite set of control signals that the quantizer can generate.
In this way, we propose to operate the controller using the cost function
\begin{align}
  \mathcal{\hat{J}} \left(\vec[][]{\beta},\vec[][]{x}(k)\right),
\end{align}
where the cost function is obtained by limiting the control signals to be part of the set $\vec{\bar{u}} \in \mathbb{W}$, where 
\begin{align}
  \mathbb{W} = \left\{\dict\vec{\beta} \, \middle| \, \vec{\beta} \in \mathcal{B} \right\}.
\end{align}
This allows us to rewrite \cref{eq:closedloopcost2ppc}, such that the quantized controller cost function becomes
\begin{align}\label{eq:closedloopcost2}
  \begin{split}
    \mathcal{\hat{J}}(\vec{\bar{u}},\vec[][]{x}(k)) &= \vec[T]{x}(k) \m[T]{\Upsilon} \m{\Upsilon} \vec{x}(k) + \vec[T]{\bar{u}}  \m{W} \vec{\bar{u}} \mathbreak 
+ 2 \vec[T]{x}(k) \m{F}\vec{\bar{u}}, ~ \vec{\bar{u}} \in \mathbb{W},
  \end{split}
\end{align}
where the minimizing $\bar{\vec{u}}$ is found by
\begin{align}\label{eq:section2new:costmin}
  \bar{\vec{u}}(\vec{x}(k)) \triangleq \argmin_{\bar{u}\in \mathbb{W}} \hat{\mathcal{J}}(\bar{\vec{u}},\vec{x}(k)).
\end{align}

The optimization problem \cref{eq:section2new:costmin} is non-convex since it operates on a finite set which is generated by the dictionary and can therefore not be minimized analytically. We therefore propose to find the minimizing solution using a greedy method, closely resembled to \ac{MP}. Thus, we select the vector from the first section in the dictionary, see \cref{fig:sparc-dictionary}, that reduces the cost function the most. This selection procedure is then repeated for every remaining section. Since $\vec[T][]{x}(k) \m[T]{\Upsilon} \m{\Upsilon} \vec[][]{x}(k) $ is constant for a given $\vec[][]{x}(k)$, this part can be neglected in \cref{eq:closedloopcost2}, reducing the cost function to
\begin{align}\label{eq:closedloopcost2reduced}
  \begin{split}
    \mathcal{\hat{J}}(\vec{\bar{u}},\vec[][]{x}(k)) &=  \vec[T]{\bar{u}} \m{W} \vec{\bar{u}}
+ 2 \vec[T]{x}(k) \m{F}\vec{\bar{u}}, \quad \vec{\bar{u}} \in \mathbb{W}.
  \end{split}
\end{align}

The greedy search is implemented for every section $m$ in the dictionary as
\begin{align}\label{eq:ppc-cost}
  &i_{m}^{*}(\vec{x}(k)) = \argmin_{i_{m}} \mathcal{\hat{J}}(\bar{\vec{u}}(i_{m}), \vec{x}(k)) \nonumber \\
  &s.t.\\
  &\bar{\vec{u}}(i_{m}) = \dictv_{i_{m}} + \sum_{j=1}^{m-1}\dictv_{i^{*}_{j}} \nonumber
\end{align}
for
\begin{align}
   &i_{m} \in ((m-1)L+1, \dots, Lm), \quad m \in \left\{1,2,\dots,M\right\} \nonumber
\end{align}
and $\dictv_{i_{m}}$ are column vectors in $\dict$. Here $L$ is the number of vectors in each section $m$ in the dictionary $\dict$. The vector $i^{*} = \left[i^{*}_{1},\, i^{*}_{2},\, \dots\, ,i^{*}_{M}\right]^{T}$ then contains the indices of the columns in the dictionary that minimize \cref{eq:closedloopcost2reduced}. This procedure is repeated for all $M$ sections in the dictionary, after which the sparse vector $\vec{\beta}(k)$ is formed as
\begin{align}
  \begin{split}
    &\vec[*]{\beta}_{j}(k) = \left\{
      \begin{array}{ll}
        1 & \text{if } j \in i^{*}\\
        0 & \text{if } j \not\in i^{*}
      \end{array}\right. ,
  \end{split}
\end{align}
such that the control signal can be reconstructed as
\begin{align}
  \vec[*]{\bar{u}}(k) = \dict \vec[*]{\beta}(k).
\end{align}

We choose to minimize the cost function for one section of $\dict$ at each iteration. The greedy search algorithm is illustrated in \cref{alg:clctrl}.

\begin{algorithm}
  \caption{Greedy search for the controller in the \ac{NCS}.}
  \label{alg:clctrl}
  \begin{algorithmic}[1]
    \State Dictionary $\dict$
    \State Input signal $\vec{x}(k)$
    \State $\vec{\beta} = \vec{0}_{ML}$
    \For{$m =1 \to M$}
    \For{$i = L(m-1)+1\to Lm$}
    \State $\bar{\vec{u}}=\dict \vec{\beta} + \dictv_{i}$
    \State $\mathit{res}_{i} = \mathcal{\hat{J}}(\bar{\vec{u}},\vec{x}(k))$
    \EndFor
    \State $g = \argmin_{i} \mathit{res}$
    \State $\vec{\beta}_{g} = 1$
    \EndFor
  \end{algorithmic}  
\end{algorithm}

\editend

The performance of the quantizer is highly dependent on the design of the dictionary, which is studied in Section \ref{sec:dictionary-design}.

Computation-wise, most of the computations in \cref{alg:clctrl} involve the evaluation of the matrix-vector and vector-vector products in the cost \cref{eq:closedloopcost2reduced}. Here the term $\vec{x}^{T}(k)\m{F}$ only has to be evaluated once at every time step. This leaves one matrix-vector and two vector-vector products that have to be evaluated $ML$ times, resulting in a total of $ML$ matrix-vector products and $2ML$ vector-vector products at every time-step. The size of the matrices $W$ and $M$ depend on the horizon length $N$. Storage-wise, the $N\times ML$ dictionary $\dict$ has to be maintained in the memory. This memory footprint can however be reduced notably since only one $N\times L$ section of the dictionary needs to be maintained in the memory for every iteration $m$. In fact one can, by generating the dictionary as described in the following sections, generate one vector of the dictionary at a time using a fixed seed for a random generator. Such a procedure would lead to a significant reduction of the memory requirements.

\begin{remark}\label{thm:section2new:greedy-alg-alternative}
  Alternatively to \cref{alg:clctrl}, one can also provide all sections at every iteration and remove the section, when one vector in the section is chosen. This can in some cases give better results, since the final cost \cref{eq:closedloopcost2reduced} might get reduced further. This though comes with a significantly increased computational cost, that increases from evaluating \cref{eq:closedloopcost2reduced} $L$ times to $L(M-m)$ times for every iteration of $m$.
\end{remark}

\subsection{Networks with two states}
\label{sec:acmjls-networks-with}

\begin{figure}[tp]
  \centering
  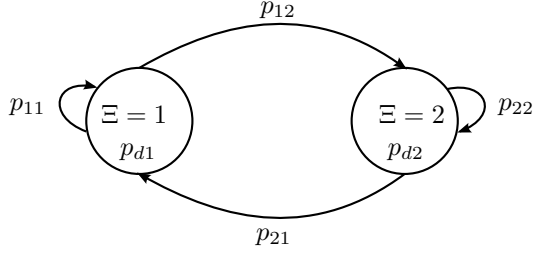
  \caption{Model of a network with two states $\netstate$, where $\netstate = 1$ corresponds to the ``good state'' featuring few dropouts, whereas $\netstate = 2$ illustrates periods with disturbances in the network, and therefore features higher packet dropout probabilities.}
  \label{fig:modeltwostate}
\end{figure}

\edit{In this section we model the network in Figure \ref{fig:clppc} by considering two network states\footnote{Extensions to multiple states present no technical difficulties.}. Here each network state represents a different probability for a packet loss to occur. This can for example be due to congestion or interference. We denote the current network state by $\netstate(k) \in \stateset$, where $\stateset = \left\{1,2\right\}$, see also \cite{pr:81,Elliott1963,6263277,DBLP:journals/corr/QuevedoJAJ13}. One state describes a reliable network situation that features low dropout probabilities, and the other state models disturbances on the network. These disturbances cause a poor connection for a period of time and therefore result in a higher dropout probability. This model is illustrated in Figure \ref{fig:modeltwostate}.}

\begin{assumption}
\label{thm:asumption1}
  The model in Figure \ref{fig:modeltwostate} is described as a discrete-time
  homogeneous Markov chain $\Xi(k)$ with transition probabilities $p_{ij} =
  \prob[\Xi(k) = i]{\Xi(k+1) = j}$ for $i,j \in \stateset$. Whether a packet
  dropout occurs or not is conditionally independent given the network state, such
  that $p_{d,i} = \prob[\Xi(k) = i]{d(k)=1}$ for $i \in \stateset$. The dropouts
  $d(k)$ are not Markovian, but correlated to $\Xi(k)$ which is
  Markovian. The augmentation of these processes $(d(k),\Xi(k))$ forms a
  Markov chain that can be classified as a \acf{MJLS}.
  
\end{assumption}

The state transition matrix for the system in Figure \ref{fig:modeltwostate} is given by
\begin{align}\label{eq:section2new:markovP}
  \m{P} = \left[
    \begin{matrix}
      p_{11} & p_{12}\\
      p_{21} & p_{22}
    \end{matrix}
    \right],
\end{align}
We assume that the current network state $\netstate(k)$ is known to both the controller and \edit{the actuator, that applies the control signal to the plant}. This can in practice be done by estimation of hidden Markov chains.


\section{Dictionary design}
\label{sec:dictionary-design}

The dictionary design has a major impact on the performance of the quantizer we proposed in Section \ref{sec:fixed-rate-quant}. A common method to design dictionaries is to use training data obtained from the signals to be quantized. This is e.g. done in the Lloyd-Max quantizer \cite{Jain:1989:FDI:59921}. An alternative method is to use the distribution of the signal to be compressed \cite{6284210,DBLP:journals/corr/abs-1212-1707}, which often is done when memoryless Gaussian sources are quantized and the distribution of the source is known or can be estimated. 

A drawback when using training data to generate the dictionary is, that the system has to operate (or is simulated) for a period of time to obtain the desired amount of training data. The simulations also have to run for a certain time to obtain training data for steady state operation. 
If only a few simulations are used, a dictionary is obtained, that only obtains information for the few special simulations. This can be omitted when averaging over multiple simulations. 
As we shall see below, when using the statistics of the \ac{NCS} and network, the dictionary can be generated offline without the need of training data. This dictionary fits the system as long as the statistics of the network do not change. \edit{The other advantage of using the system statistics for the design of the dictionary is that the dictionary can be designed alongside the controller, and easily can be modified when the controller is tuned.}

We propose to use the distribution of the \ac{NCS} and the statistics of the network to generate the dictionary offline. Using this, we can design a dictionary that fits the system when the first and second moment of the system are known or can be estimated. 

\edit{For the ease of exposition we first describe how to design the dictionary for a single network state. This will in \cref{sec:dict-design-two} be expanded to take the state transition probabilities in $P$ into account to design a single dictionary that covers multiple network states.}

\subsection{Dictionaries considering a single network state}
\label{sec:dict-single-state}

In this section, we describe the design of the dictionary for a \ac{NCS} with a single network state as described in Section \ref{sec:acppc}. The probability for a packet dropout to occur is $p_{d}$.

For the setup considered, we note that the first and second moment of the state $\vec[][]{x}(k)$ of the \ac{NCS} are given by

\begin{align}
&\expval{\vec[][]{x}(k)} = \expval{\vec[][]{x}(k) } =  \expval{ \vec[][]{x}(k)} = \vec{0}_{p}\label{eq:expQu}\\
  &\m[][\bm{x}]{Q} = \operatorname{var}\left\{ \vec[][]{x}(k) \right\} =  \expval{\vec[][]{x}(k)\vec[T][]{x}(k)}. \label{eq:varQu}
\end{align}

Using (\ref{eq:nextstate}) and (\ref{eq:buffer}), we describe the aggregated system state 
\begin{align}\label{eq:theta}
  \vec[][]{\Theta}(k) = \left[
    \begin{matrix}
      \vec[][]{x}(k)\\
      \vec[][]{b}(k-1)
    \end{matrix}
\right] \in \R^{p+N}
\end{align}
at time instance $k$ by
\begin{align} \label{eq:sysstability}
  \vec[][]{\Theta}(k+1) = \bar{\m{A}}(d(k))\vec[][]{\Theta}(k) + \bar{\vec[][]{B}}\omega(k)
\end{align}
where
\[
\begin{minipage}{0.48\linewidth}
  \raggedleft
\begin{align*}
  \m[][]{\bar{A}}(0) = \left[
    \begin{matrix}
      \m{A} - \vec[][1]{B}\vec[T][1]{e}\m{K} & \bm{0}_{p\times N}\\
      -\m{K} & \bm{0}_{N\times N}
    \end{matrix}
  \right]
\end{align*}
\end{minipage}
\begin{minipage}{.48\linewidth}
\begin{align*}
  \m[][]{\bar{A}}(1) = \left[
    \begin{matrix}
      \m{A} & \vec[][1]{B}\vec[T][1]{e}\m{M} \nonumber \\
     \bm{0}_{N\times p} & \m{M} 
    \end{matrix}
  \right]
\end{align*}
\end{minipage}
\]
\begin{align*}
  \vec[][]{\bar{B}}(0) =   \vec[][]{\bar{B}}(1)  = \left[
    \begin{matrix}
      \vec[][2]{B}  \\ 
      \bm{0}_{N} 
    \end{matrix}
  \right] = \vec{\bar{B}},
\end{align*}
and $d_{k} = 1$ when a packet dropout occurred. 

\begin{definition}{\acs{MSS}\cite{5742978,DT-MJLS}}
\edit{The linear system \cref{eq:sysstability} is \ac{MSS} if there exist a bounded $\mu$ and $\m{Q}_{\m{\Theta}}$, such that 
  \begin{align}
    &\expval{\vec{\Theta}(k)} \rightarrow \mu, &\quad k \rightarrow \infty\\
    &\expval{\vec{\Theta}(k)\vec{\Theta}^{T}(k)} \rightarrow \m{Q}_{\m{\Theta}}, &\quad k \rightarrow \infty \label{eq:section3:MSS-def-var}
  \end{align}
for all initial conditions $\vec{\Theta}(0)$ with bounded variance and $d(0)\in\left\{0,1\right\}$.
}
\end{definition}

\edit{Thus for system \cref{eq:sysstability} to be \ac{MSS}, we require that the first and second moments of \cref{eq:sysstability} converge to a finite value. When the packet dropouts are \ac{IID}, the recursion \cref{eq:sysstability} and the distribution of $d(k)$ amount to a \ac{MJLS} with transition matrix \cite{5742978} 
  \begin{align}
\left[  \begin{matrix}
    1-p_{d} & p_{d}\\
    1-p_{d} & p_{d}
  \end{matrix}\right]\label{eq:section3:1-p_d-}.
  \end{align}
As described in \cite[Proposition 3.6]{DT-MJLS}, the \ac{MJLS} is stable in its first moment if it is stable in its second moment. We therefore only have to verify that the second moment of \cref{eq:sysstability} converges. Since $\omega(k)$ is \ac{IID} and thereby wide sense stationary, it is sufficient to verify that the homogeneous system $\m{\Theta}(k+1)= \m{\bar{A}}(d(k))\m{\Theta}(k)$ is \ac{MSS} \cite[Theorem 3.33]{DT-MJLS}. The second moment of \cref{eq:sysstability} at time $k+1$ is given by
  \begin{align}
    \begin{split}\label{eq:section3:MSS-expThetaCalc}
      \m{Q}_{\m{\Theta}}(k+1) &= \expval{\m{\Theta}(k+1)\m{\Theta}^{T}(k+1)} \\
      &= \expval{\m{A}(d(k)) \m{Q}_{\m{\Theta}}(k) \m{A}^{T}(d(k))},
    \end{split}
  \end{align}
which can be rewritten using properties of the Kronecker product \cite{DT-MJLS} to
\begin{align}
  \operatorname{vec}\left\{\m{Q}_{\m{\Theta}}(k+1)\right\} = \m{\Psi}\operatorname{vec}\left\{\m{Q}_{\m{\Theta}}(k)\right\} \in \R^{(p+N)^{2}},
\end{align}
where $\operatorname{vec}\left\{\m{Q}\right\}$ stacks the columns of $\m{Q}$ and
\begin{align}
  \begin{split}\label{eq:section3:secondmoment-stab}
    \m{\Psi} &= \expval{\m{A}(d(k))\kron \m{A}(d(k))} \in \R^{(p+N)^{2}\times (p+N)^{2}} \\
    &= \left(1-p_{d}\right) \bar{\m{A}}(0) \kron \bar{\m{A}}(0) + p_{d} \bar{\m{A}}(1) \kron \bar{\m{A}}(1),
  \end{split}
\end{align}
with $\kron$ being the Kronecker product. The system \cref{eq:sysstability} is \ac{MSS} if  $\m{Q}_{\m{\Theta}}(k)\rightarrow 0$ as $k\rightarrow \infty$. This occurs only if $p_{d}$ is such that all eigenvalues of $\m{\Psi}$ are within the unit circle.
}

\edit{In this work we only consider \ac{MSS} systems, which requires that there exists a stabilizing feedback gain \cref{eq:section2new:ctrlGain} and that $p_{d}$ is such that \cref{eq:section3:secondmoment-stab} has its eigenvalues within the unit circle. Note that \ac{MSS} can be shown using other methods as described in \cite[Theorems 3.9 and 3.33]{DT-MJLS} and \cite[Theorem 2]{5742978}.}

\begin{lemma}
  \label{thm:1}
  If the \ac{NCS} in (\ref{eq:sysstability}) is \ac{MSS}, the variance of $\omega(k)$ is finite and, in the absence of quantization effects, the first- and second-order moments of the state in (\ref{eq:sysstability}) are given by (\ref{eq:expQu}) and (\ref{eq:varQu}), then 
\begin{align}
  \m[][\bm{\Theta}]{Q} = \lim_{k\rightarrow
    \infty}\expval{\m[][]{\Theta}(k)\m[T][]{\Theta}(k)},\label{eq.mbmth-=-lim_kr}
\end{align}
which can be computed as
\begin{align}\label{eq:thetarecursive}
  \m[][\bm{\Theta}]{Q} = \m{\mathcal{A}}\m[][\bm{\Theta}]{Q}\m[T][]{\mathcal{A}} + p_{d}(1-p_{d})\m{\tilde{\mathcal{A}}}\m[][\bm{\Theta}]{Q}\m[T]{\tilde{\mathcal{A}}} + \sigma_{\omega}^{2}\vec{\bar{B}}\vec[T]{\bar{B}},
\end{align}
where
\begin{align}\label{eq:section3:expvalmarkovA}
    \m{\mathcal{A}} = \expval{\m{\bar{A}}(d(k))} = p_{d}\bar{\m{A}}(1) + (1-p_{d})\bar{\m{A}}(0)
\end{align}
and
\begin{align}
  \m{\mathcal{\tilde{A}}} = \m{\bar{A}}(1) - \m{\bar{A}}(0).
\end{align}

\end{lemma}

\noindent\textit{Proof:} see Appendix \ref{sec:proof-lemma-}.

\edit{The second moment of the aggregated system state $\vec{\Theta}(k)$ found in \cref{thm:1} can be used to design the vectors in the dictionary $\dict$ using a Gaussian distribution, such that
\begin{align}
  \dictv_{i} \sim \mathcal{N}\left(\bm{0}_{N},\m[][\bm{u}]{Q}\right), \quad i\in\left\{1,2,\dots,ML\right\} \label{eq:distdict}
\end{align}
where 
\begin{align}\label{eq:covucalc}
  \m[][\bm{u}]{Q} = \m{K}
  \left[\begin{matrix}
    \m[][p]{I} & \m{0}_{N}^{T}
  \end{matrix}\right]
\m[][\bm{\Theta}]{Q}
  \left[\begin{matrix}
    \m[][p]{I} \\
    \m{0}_{N}
  \end{matrix}\right]
 \m[T]{K}.
\end{align}
}

The above result gives statistics of the \ac{NCS} operating while neglecting quantization effects. By creating the dictionary using (\ref{eq:distdict}) and using a bit rate that is high enough, the granular distortion added by the quantizer is small compared to $\sigma_{\omega}^{2}$ and will only have a limited impact on the performance of the \ac{NCS}. When the bit rate is decreased, we can compensate for the granular distortion by scaling the dictionary with a factor larger than 1.

We design two dictionaries using (\ref{eq:distdict}). The first dictionary is a \ac{GR} dictionary where every section $m\in \left\{0,1,\dots,M-1\right\}$ is generated using (\ref{eq:distdict}). In the second dictionary (the \ac{GSR} dictionary), we, inspired by \cite{DBLP:journals/corr/abs-1212-1707}, scale every section $m\in \left\{0,1,\dots,M-1\right\}$. 
\edit{In \cite{DBLP:journals/corr/abs-1212-1707}, the asymptotically optimal scaling factor was found to be
  \begin{align}
    c_{m} = a_{0}^{\frac{1}{2}} a_{1}^{\frac{m}{2}}
  \end{align}
for some $a_{0}$ and $a_{1}$ that are independent of $m$. However, in non-asymptotical cases, we have experimentally observed that a better choice is to use 
\begin{align}\label{eq:dictscaling}
  c_{m} = a_{0}'^{\frac{1}{2}} a_{1}'^{\frac{m}{M}},
\end{align}
where $a_{0}'= 1$ and $a_{1}' = \frac{1}{M}$. The scaling factor provides a trade-off between granular and overload distortions. Since the variance of the signal to be quantized is decreased in each iteration (as $m$ increases), the scaling factor ensures that the variance of the dictionary elements is decreased accordingly to better match the changing input statistics.
}

The vectors in the \ac{GSR} dictionary can then be created using the scaling factor $c_{m}$ from \cref{eq:dictscaling} and a Gaussian distribution, such that
\begin{align} \label{eq:gsr-dict}
    \dictv_{mL+i} \sim c_{m}\mathcal{N}\left(\m{0}_{N},\m[][\bm{u}]{Q}\right),
\quad 
  \begin{split}
    &i\in \left\{1,2,\dots,L\right\}, \\
    &m \in \{0,1,\dots,M-1\}.
  \end{split}
\end{align}
\edit{Using this scaling factor, the variance of the individual sections in $\dict$ decreases as $m$ increases. The idea is, that using the greedy algorithm (\cref{alg:clctrl}) on the \ac{GSR} dictionary, the vector that is selected first will reduce the residual in \cref{eq:section2new:quant-error} the most. The vector that is selected hereafter will reduce it slightly less and so on. This is expected to give a lower residual than the \ac{GR} dictionary, where every section in $\dict$ has the same variance.}

\begin{remark}\label{thm:section3:editcrefthm:2-dictionaries}
  \edit{When there are multiple dropout scenarios, as presented in \cref{sec:acmjls-networks-with}, the above described method can be used to design a dictionary for each network state. Using this approach, the controller and actuator switch dictionaries when the network changes from one state to the other. This requires that both the controller and actuator know the current network state. The dictionaries can in this case be designed using \cref{eq:distdict,eq:gsr-dict} where the covariance is calculated as in \cref{thm:1}. Here the first dictionary is designed using $p_{d} = p_{d1}$ and the second is designed using $p_{d}=p_{d2}$. In the implementation of this, \cref{alg:clctrl2} is used to switch between the dictionaries in the \ac{NCS}. Note however, that although the dictionaries designed using this method take the individual dropout probabilities in the current network state into account, they do not consider the transition probabilities to go from one state to the other.}
\end{remark}

\begin{algorithm}
  \caption{Modification to \cref{alg:clctrl} to perform greedy search with two network states.}
  \label{alg:clctrl2}
  \begin{algorithmic}[1]
    \State Dictionaries $\dict_{1}, \dict_{2}$
    \State Input signals  $\vec{x}(k), \netstate(k)$ 
    \State $\dict = \dict_{\netstate(k)}$
    \State Go to step 3 of \cref{alg:clctrl}.
  \end{algorithmic}
\end{algorithm}


\subsection{Dictionary design for two-state networks}
\label{sec:dict-design-two}

We will in this section describe a single dictionary that can be used when there are multiple network states as described in \cref{sec:acmjls-networks-with}. The main advantage here compared to the \edit{method described in \cref{thm:section3:editcrefthm:2-dictionaries} is}, that the controller and buffer do not need to have any information on the current network state. This dictionary uses the statistics of the system featuring the network model described in Section \ref{sec:acmjls-networks-with}. The second moment is found by adopting the methods described in \cite{6263277} and \cite{DBLP:journals/corr/QuevedoJAJ13}. The quantizer noise is not taken into account in this analysis, and it therefore only shows the performance of the \ac{NCS} without any quantization. This way, under Assumption \ref{thm:asumption1}, the \ac{NCS} shown in Figure \ref{fig:clppc} can be described by the jump-linear model (\ref{eq:sysstability}) based on (\ref{eq:nextstate}) and (\ref{eq:buffer}) where $d(k) = 1$ indicates that a packet dropout occurred at time instance $k$, and $\Theta(k+1)$ is described by \edit{the recursion} (\ref{eq:sysstability}).
These models have been studied in e.g. \cite{DT-MJLS}, where Theorems 3.9 and 3.33 from \cite{DT-MJLS} describe necessary and sufficient conditions for \ac{MSS}.

\begin{lemma}\label{thm:2}
  If the system in (\ref{eq:sysstability}) is \ac{MSS}, Assumption
  \ref{thm:asumption1} holds, quantizer effects are neglected, the variance of $\omega(k)$ is finite and
the second moment of the state in \cref{eq:sysstability} is given by (\ref{eq:varQu}), then
  \begin{align}
    \label{eq:cov2state}
    \m[][\bm{\Theta}]{Q} = \lim_{k \to \infty} \expval{\m[][]{\Theta}(k) \m[T][]{\Theta}(k)} =
    \sum_{j \in \stateset} \m[][j]{Q} 
  \end{align}
where
\begin{align}
  \label{eq:6}
  \m[][j]{Q} =  \sum_{i \in \stateset} p_{ij} \m[ ][j]{\mathcal{A}} \m[ ][i]{Q} \m[T][j]{\mathcal{A}} + \pi_{j}\sigma_{\omega}^{2} \vec{\bar{B}}\vec[T]{\bar{B}},
\end{align}
with $\pi_{j}$ being the stationary probability of the Markov state $\netstate = j$ and
\begin{align}\label{eq:dictionary2:markovA}
  \m[][j]{\mathcal{A}} &= \expval[\Xi(k) = j]{\m{\bar{A}}(d(k))} \nonumber \\
  & = p_{d,j}\bar{\m{A}}(1) + (1-p_{d,j})\bar{\m{A}}(0) \hspace{0.5cm} j \in \stateset.
\end{align}
%
%
%
%
\end{lemma}

\noindent\textit{Proof:} see Appendix \ref{sec:proof-theorem2}.

The dictionaries for the \ac{NCS} featuring correlated dropouts can thus be described using \cref{eq:cov2state}. Using this, we design the \ac{GR2} dictionary using \cref{eq:distdict,eq:covucalc} with covariance (\ref{eq:cov2state}) and the \ac{GSR2} dictionary where every section $m \in M$ is scaled using \cref{eq:gsr-dict}. \edit{This dictionary can be implemented in the \ac{NCS} using \cref{alg:clctrl}. It should further be noted, that this design does not require that the controller and actuator have knowledge of the current network state, merely statistical information is used.}



\section{Simulations}
\label{sec:simulations}
\acused{MSE}

In this section we provide simulation studies for the fixed rate quantizer on \acp{NCS}. We design the dictionaries as described in Section \ref{sec:dictionary-design} and compare these to the \ac{NCS} running without quantizer and the \ac{NCS} with the quantizer using a dictionary containing Gaussian \ac{IID} samples. The latter is a simple way to generate dictionaries when quantizing Gaussian \ac{IID} sources \cite{6284210,DBLP:journals/corr/abs-1212-1707} and is expected not to achieve a competitive performance compared to the \ac{GR} and \ac{GSR} dictionaries.  
\edit{We compare the bit rate in bit/symbol against the empirical cost of the state in the system and the control signal. This cost is calculated as 
\begin{align}
  \mathrm{MSE} = \|\vec{x}(k)\|_{\m{Q}}^{2} + \|\vec{u}(k)\|_{R}^{2},
\end{align}
which is averaged over the length of the simulation.}
When the bit rate is too low, the quantizer may overload and make the system unstable. 

\edit{In the simulations we use the recursion \cref{eq:nextstate} where the noise $\omega(k) \sim \mathcal{N}(0,\sigma_{\omega}^{2})$ is zero mean Gaussian with variance $\sigma_{\omega}^{2} = 1$, and} the system matrix $\m{A}$ is randomly generated as
\begin{align}\label{eq:simA}
  \m{A} = \left[
    \begin{matrix}  
  -0.758&-0.325&-0.085&0.060&-2.256\\
  0.432&-0.356&0.002&0.007&-0.171\\
  -0.173&1.063&0.366&0.671&0.939\\
  0.951&0.667&0.737&-0.434&0.352\\
  1.054&0.484&-0.158&0.454&-0.264
    \end{matrix} \right].
\end{align}
\edit{$\m{A}$ has absolute values of the eigenvalues $1.659, \,1.659 ,\, 1.241, \,0.754$ and $0.43$}. The matrices $\m[][1]{B} = \m[][2]{B} = \left[1\,1\, 1\, 1\, 1\right]^{\mathrm{T}}$, \edit{such that the noise and the control signal affect all states}. \edit{This system is fully controllable. }\edit{The weighting parameters used in \cref{eq:section2new:weights}, that are used in the cost function \cref{eq:closedloopcost2reduced} are $\m{Q} = \m{I}_{5}$ and $R = 1$}. The horizon length $N$ is set to 5 and the \edit{dictionary $\dict$, shown in \cref{fig:sparc-dictionary},} contains $M=2$ sections, since this shows the best performance in simulations and maintains stability at lower bit rates than a dictionary containing 3 or more sections. The number of vectors $L$ in each section is calculated from the desired bit rate $R$ by
\begin{align}
  \label{eq:calcL}
  L = \ceil{2^{\frac{N R}{M}}},
\end{align}
which is rounded towards the greatest integer. The actual bit rate is hereafter calculated by solving (\ref{eq:calcL}) with respect to $R$. \edit{The dictionary $\dict$ is then designed as described in \cref{thm:1}, where the \ac{GR} dictionary is generated using \cref{eq:distdict} and the \ac{GSR} dictionary using \cref{eq:gsr-dict}.}  
\edit{The control signal $\vec{u}(k)$ is then calculated using \cref{alg:clctrl}}

We run multiple simulations each of \SI{50000}{} time instances over which we average the \ac{MSE}. The dictionaries are randomly generated and the dropouts are randomized for each simulation. If one simulation with a certain dictionary is unstable at a given bit rate, all simulations for this dictionary at this bit rate are considered to be unstable.

The dictionaries are scaled by a fixed scaling that results in the best performance. The \ac{IID} dictionary is generated using a Gaussian distribution with $\sigma^{2}=25$, which provides the best performance in \ac{MSE} compared to bit rate without being overloaded at lower bit rates. The \ac{GSR} dictionary is scaled by a factor 2, whereas the \ac{GR} dictionary is scaled by 1. 

\subsection{\ac{NCS} with \ac{IID} dropouts}
\label{sec:acncs-with-aciid}

\begin{figure}[tp]
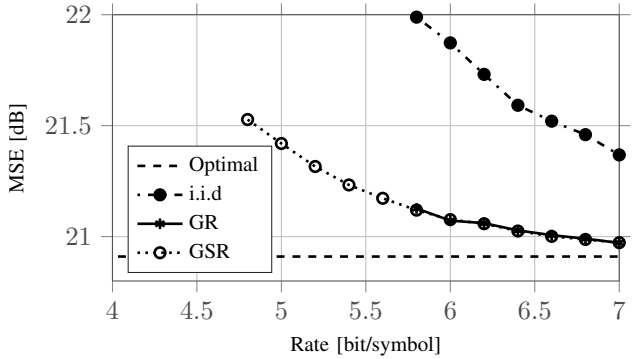

  \centering
  \inputtikz{simulations/fig/mse_rate_singledrop_singlesys2/}{simsingle0.1}
  \caption{Average \ac{MSE} over 12 simulations. The designed dictionaries with stars showing the \ac{GR} dictionary and circles the \ac{GSR} dictionary. These are compared to the \edit{optimal performance} of the \ac{NCS} without quantizer (dashed) and an \ac{IID} (dash-dot) dictionary.}
  \label{fig:single11avg}
\end{figure}

We first simulate a \ac{NCS} with a single network state featuring \ac{IID} random packet dropouts with the probability $p_{d}$, which is set to $0.10$ in this simulation. Figure \ref{fig:single11avg} shows the average \ac{MSE} over 12 simulation runs. 

The results show that the \ac{GSR} dictionary (that is generated using (\ref{eq:gsr-dict})) can maintain stability at bit rates down to \SI{4.8}{\bit\per\symbol}. The \ac{MSE} is less than \SI{1}{\decibel} higher than the \ac{NCS} without any type of quantizer, and is reduced to only \SI{0.1}{\decibel} at \SI{7}{\bit\per\symbol}. 
This is a significant performance improvement compared to the \ac{IID} Gaussian dictionary where the \ac{MSE} reaches \SI{21.6}{\decibel} already at \SI{6.4}{\bit\per\symbol}, while the \ac{GSR} dictionary stays below this even at a rate of \SI{4.8}{\bit\per\symbol}.
The \ac{GR} dictionary (generated using (\ref{eq:distdict})) maintains stability at bit rates of \SI{5.8}{\bit\per\symbol}, where it results in identical performance, with respect to \ac{MSE} as the \ac{GSR} dictionary. 

\begin{figure}[t]
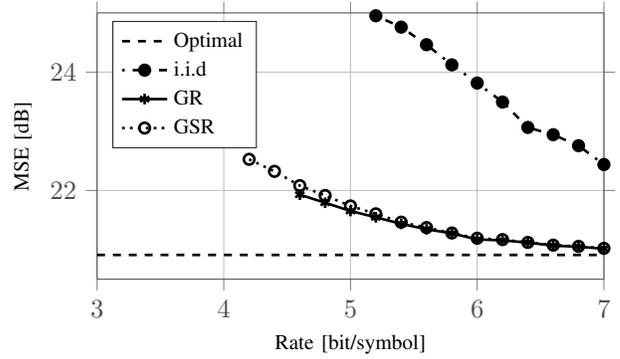

  \centering
  \inputtikz{simulations/fig/mse_rate_singledrop_singlesys2/}{simsingle0.2}
  \caption{Average \ac{MSE} over 12 simulations. The designed dictionaries with stars showing the \ac{GR} dictionary and circles the \ac{GSR} dictionary. These are compared to the \edit{optimal performance} of the \ac{NCS} without quantizer and an \ac{IID} dictionary.}
  \label{fig:single12avg}
\end{figure}

Figure \ref{fig:single12avg} shows the results when the dictionaries use an additional scaling factor of 2, such that the \ac{GR} dictionary is scaled by a factor 2 and the \ac{GSR} is scaled by a factor 4. The \ac{MSE} is slightly increased compared to Figure \ref{fig:single11avg} for higher bit rates, but the system is stable at lower bit rates. The \ac{GSR} dictionary is able to maintain stability at bit rates down to \SI{4.2}{\bit\per\symbol}. The \ac{GR} dictionary results in a slightly lower \ac{MSE}, but is unable to maintain stability at bit rates lower than \SI{4.6}{\bit\per\symbol}. The \ac{MSE} of the \ac{IID} Gaussian dictionary increased compared to Figure \ref{fig:single11avg}, but is able maintain stability at bit rates down to \SI{5.2}{\bit\per\symbol}.

\edit{Both simulations show that all dictionaries asymptotically approach the optimal \ac{MSE} of the unquantized \ac{NCS} as the bit rate increases. In \cref{fig:single12avg}, there is only \SI{0.11}{\decibel} difference in the \ac{MSE} between the optimal \ac{MSE} from the \ac{NCS} with no quantizer and the \ac{GR} and \ac{GSR} dictionaries. 
}

\subsection{\ac{NCS} with correlated dropouts}
\label{sec:acncs-with-two}

We simulate a system featuring a network that is modeled as shown in Figure \ref{fig:modeltwostate} using two network states. The transition matrix is defined by
\begin{align}
  \label{eq:transm}
  \m{P} = 
    \begin{bmatrix}
      0.95 & 0.05\\
      0.25 & 0.75
    \end{bmatrix},
\end{align}
and the dropout probabilities $p_{d1}= 0.05$, $p_{d2}= 0.15$. \edit{Here the dictionary design for the \ac{GR} and \ac{GSR} dictionaries is done as described in \cref{thm:section3:editcrefthm:2-dictionaries}. Using these dictionaries, the control signal $\vec{u}(k)$ is computed using \cref{alg:clctrl2}.}
\edit{For the \ac{GR2} and \ac{GSR2} dictionaries, the dictionary design is done as described in \cref{thm:2}, where the \ac{GR2} dictionary then is generated using \cref{eq:distdict} and  the \ac{GSR2} dictionary using \cref{eq:gsr-dict}. The control signal $\vec{u}(k)$ is then computed using \cref{alg:clctrl}.}

\begin{figure}[t]
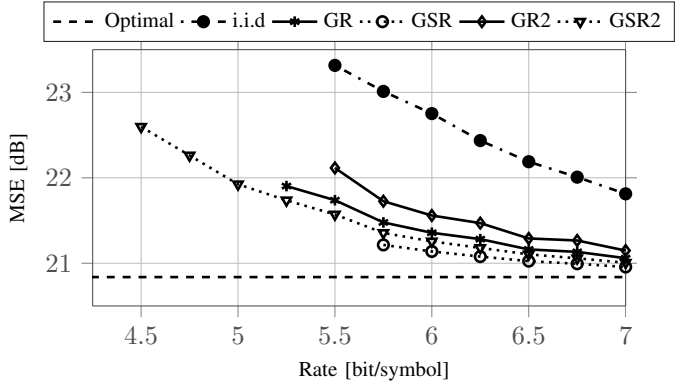

  \centering
  \inputtikz{simulations/fig/mse_rate_twostate_new/}{twostate}
  \caption{Average \ac{MSE} over 24 simulations. The designed dictionaries with stars showing the \ac{GR} dictionary and circles the \ac{GSR} dictionary. The diamonds and triangles show the \ac{GR2} and \ac{GSR2} dictionary, respectively. These are compared to the \edit{optimal performance} of the \ac{NCS} without quantizer and an \ac{IID} dictionary.}
  \label{fig:twostate12avg}
\end{figure}

The dictionaries overload when these are scaled with the earlier mentioned factors, and are thereby unable to maintain stability. 
We therefore change the scaling factors such that the stability and the \ac{MSE} at a rate of \SI{7}{\bit\per\symbol} are maintained for the different dictionaries.
Thus, the \ac{GR} dictionary is scaled by 2, the \ac{GR2} and \ac{GSR2} dictionaries are scaled by 3, the \ac{GSR} dictionary is scaled by 4, and the \ac{IID} Gaussian dictionary is scaled by 2, such that here $\sigma^{2}=100$. 
Figure \ref{fig:twostate12avg} shows the results averaged over 24 simulations. 

The \ac{GR} dictionary shows a better performance at low bit rates compared to he \ac{GSR} dictionary, but generally has a higher \ac{MSE}. The \ac{GSR2} dictionary shows the far best performance, being able to maintain the system stable at bit rates down to \SI{4.5}{\bit\per\symbol}, which is far below any other dictionary. The \ac{GR2} dictionary performs slightly worse than the \ac{GR} dictionary. It should here though be noted, that the \ac{GR} and \ac{GSR} dictionaries assume that the controller and buffer know the network state at every time instance, whereas the \ac{GR2} and \ac{GSR2} do not need this information. The dictionary containing \ac{IID} samples generally has a higher \ac{MSE} at comparable bit rates to the designed dictionaries.

The simulations generally show, that the risk of overloading the quantizer is reduced when the scaling of the dictionaries is increased. This scaling though results in a higher \ac{MSE} for comparable bit rates. 

\edit{The simulations show that the proposed dictionaries result in a significantly improved performance compared to a dictionary consisting of simple \ac{IID} Gaussian generated variables. The performance of the quantized \ac{NCS} is however highly dependent on the system matrix $\m{A}$, the variance of the system disturbance $\sigma_{\omega}^{2}$ and the packet dropout probability $p_{d}$. The system matrix $\m{A}$ and $p_{d}$ are linked through the second moment of the recursion \cref{eq:sysstability} through \cref{eq:section3:secondmoment-stab}. 
Here it is important to mention that the system matrix $\m{A}$ is determined by the dynamics of the system to be controlled \cite{astrom:2007:feedbacksystems}. Further, the packet dropout probability $p_{d}$ depends on the network (network load, interference etc.). The higher the eigenvalues in \cref{eq:section3:secondmoment-stab} the more ``flat'' the distribution \cref{eq.mbmth-=-lim_kr} becomes. This results in more granular distortion since the entries in the dictionary are more wide-spread, which then results in a higher \ac{MSE} and a higher bit rate is required to maintain stability.}

\begin{figure}[tp]
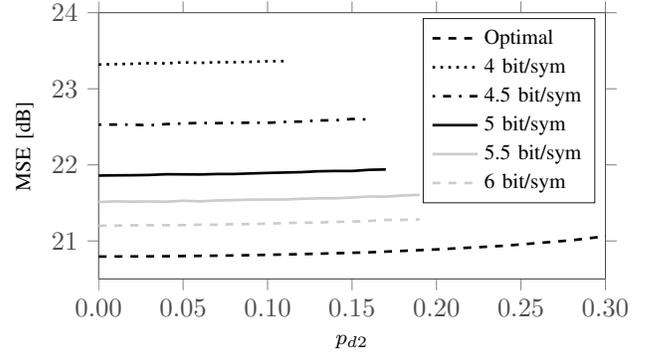

  \centering
  \inputtikz{simulations/fig/sweeppd2/}{sweeppd2}
  \caption{\edit{The \ac{MSE} shown at different rates for the \ac{GSR2} dictionary and the system using the optimal control input. Here $p_{d2}$ is swept while $p_{d1}$ and $P$ are kept constant. }}
  \label{fig:simulations:cent-capt-swept}
\end{figure}

\edit{The impact on the performance of the quantized \ac{NCS} when the network is considered Markovian is further affected by the transition matrix $\m{P}$, which also depends on the network through which the controller and system are connected. The impact of the transition probabilities in $\m{P}$ on the performance of the quantized \ac{NCS} directly depends on the underlying packet dropout probabilities $p_{d1}$ and $p_{d2}$. To illustrate the impact of the packet dropout probabilities affecting the ``bad'' state of the network, we sweep $p_{d2}$ in simulations while $\m{P}$ is maintained as in \cref{eq:transm}, $p_{d1} = 0.5$ and $\m{A}$ is as in \cref{eq:simA}. The results of this are shown in \cref{fig:simulations:cent-capt-swept}. The figure illustrates that stability can be maintained at higher packet dropout rates in the ``bad'' network state as the bit rate increases.  The reason for this is that the system is less affected by the distortion that is introduced when the quantizer is overloaded. }

\section{Conclusions}
\label{sec:conclusions}

We have presented a quantized controller for \ac{NCS}, that features a fixed rate \ac{VQ}. The performance of this setup is highly dependent on the design of the dictionaries used. Simulations show, that the proposed dictionaries perform significantly better than a dictionary containing \ac{IID} random samples. Performance approaches the \ac{MSE} of the \ac{NCS} without quantizer when the bit rate increases. The \ac{GSR} dictionary generally shows a slightly better performance and overloads at lower bit rates than the \ac{GR} dictionary. 
When assuming correlated dropouts, the \ac{GR} and \ac{GSR} dictionary also outperform the \ac{IID} dictionary. This scenario is also tested using the \ac{GR2} and \ac{GSR2} dictionaries, that utilize a \ac{MJLS} model to describe the stationary covariance of the \ac{NCS}, where the \ac{GSR2} performance generally outperforms the other dictionaries. An additional advantage with the \ac{GR2} and \ac{GSR2} dictionaries is, that the controller and buffer do not need to know in which state the network is. 


\appendices
\setcounter{equation}{0}
\renewcommand{\theequation}{\Alph{section}.\arabic{equation}}

\section{Proof of Lemma \ref{thm:1}}
\label{sec:proof-lemma-}

Proceeding as in \cite{6416000}, the covariance of the \ac{PPC} is based on the dropout probability $p_{d}$, which in this case is known. 

The covariance of $\m[][]{\Theta}(k)$ in (\ref{eq:sysstability}) can then be described by
\begin{align}\label{eq:covarianceSection}
  \m[][\Theta]{Q}(k+1) &= \mathbb{E}\left\{\m[][]{\Theta}(k+1)\m[T][]{\Theta}(k+1)\right\} \nonumber \\
&= \m{\mathcal{A}}\mathbb{E}\left\{\m[][]{\Theta}(k)\m[T][]{\Theta}(k)\right\} \m[T]{\mathcal{A}} +
 \\
& ~~~ p_{d}(1-p_{d})\m{\tilde{\mathcal{A}}}\mathbb{E}\left\{\m[][]{\Theta}(k)\m[T][]{\Theta}(k)\right\}\m[T]{\tilde{\mathcal{A}}} + \sigma_{w}^{2}\vec[][w]{\mathcal{B}}\vec[T][w]{\mathcal{B}},\nonumber
\end{align}
where
\begin{align}
  &\m{\mathcal{A}} = \mathbb{E}\left\{ \m{\bar{A}}(d(k)) \right\} = p_{d}\bar{\m{A}}(1) + (1-p_{d})\bar{\m{A}}(0)\\
  &\m{\mathcal{\tilde{A}}} = \m{\bar{A}}(1) - \m{\bar{A}}(0)\\
  &\m{\mathcal{B}} = \mathbb{E}\left\{\vec{\bar{B}}(d(k)) \right\} = \bar{\vec{B}}.
\end{align}
We rewrite (\ref{eq:covarianceSection}) to
\begin{align}\label{eq:appcov}
  \begin{split}
    &\m[][\Theta]{Q}(k+1) = \m{\mathcal{A}}\m[][\Theta]{Q}(k) \m[T]{\mathcal{A}} \mathbreak
+ p_{d}(1-p_{d}) \m{\tilde{\mathcal{A}}} \m[][\Theta]{Q}(k) \m[T]{\tilde{\mathcal{A}}}  + \sigma_{w}^{2}\vec[][w]{\mathcal{B}}\vec[T][w]{\mathcal{B}},
  \end{split}
\end{align}
where the stationary covariance can be found by
\begin{align}
  \m[][\Theta]{Q} = \lim_{k\rightarrow \infty}\m[][\Theta]{Q}(k).
\end{align}
This can be solved by finding the symmetric $\m[][\Theta]{Q} > 0$ that satisfies the linear matrix equation (\ref{eq:thetarecursive}), which either can be done by iterating (\ref{eq:appcov}) or by using the closed form solution described in Remark 2 in \cite{6416000}.

The matrix $\m[][\Theta]{Q}$ contains the covariances of the the state $\vec{x}$ and buffer $\vec{b}$ from (\ref{eq:theta})
\begin{align}\label{eq:Qtheta}
  \m[][\bm{\Theta}]{Q} = \left[
    \begin{matrix}
      \m[][\bm{x}]{Q} & \mathbb{E}\left\{\vec{x}\vec[T]{b}\right\}\\
      \mathbb{E}\left\{\vec{b}\vec[T]{x}\right\} & \m[][\bm{b}]{Q}
    \end{matrix}
\right].
\end{align}
We are interested in the covariance of $\vec{u}$, which according to (\ref{eq:varQu}) depends on the covariance of $\vec[][]{x}(k)$, denoted $\m[][x]{Q}$, which is the upper diagonal part of $\m[][\Theta]{Q}$

By isolating the covariance $\m[][x]{Q}$ of the plant state $\vec{x}$, the covariance of the controller output $\vec{u}$ is given by
\begin{align}\label{eq:covucalcapp}
  \m[][\bm{u}]{Q} = \m{K}
  \left[\begin{matrix}
    \m[][p]{I} & \m{0}
  \end{matrix}\right]
\m[][\bm{\Theta}]{Q}
  \left[\begin{matrix}
    \m[][p]{I} \\
    \m{0}
  \end{matrix}\right]
 \m[T]{K}.
\end{align}


\section{Proof of Lemma \ref{thm:2}}
\label{sec:proof-theorem2}

This proof follows the procedure of \cite{DBLP:journals/corr/QuevedoJAJ13}. Using the law of total expectation, we have
\begin{align}
  \expval{\m[][]{\Theta}(k+1)\m[T][]{\Theta}(k+1)} = \sum_{j\in \mathbb{B}} \m[][j]{Q}(k+1), \label{eq:covQ}
\end{align}
with 
\begin{align}\label{eq:BQ}
  \begin{split}
    &\m[][j]{Q}(k+1) = \expval[\netstate(k)=j]{ \m[][]{\Theta}(k+1) \m[T][]{\Theta}(k+1)} \mathbreak
\times \prob{\netstate(k)=j}
  \end{split}
\end{align}
Using assumption \ref{thm:asumption1} on (\ref{eq:sysstability}), we write
\begin{align}
  &\expval[\netstate(k)=j]{\m[][]{\Theta}(k+1)\m[T][]{\Theta}(k+1)}  \nonumber \\
  & ~= \expvalSym \left\{ \left(\bar{\m{A}}(d(k))\m[][]{\Theta}(k) + \bar{\m{B}}\omega(k) \right)\left(\bar{\m{A}}(d(k))\m[][]{\Theta}(k) \right.\right.\\
&\qquad \left. \left. + \bar{\m{B}}\omega(k) \right)^{\mathrm{T}}\middle| \netstate(k)=j \right\}, \nonumber
\end{align}
which can, considering $\omega(k)$ being white Gaussian noise, be rewritten to
\begin{align}
  \begin{split}
    &\expval[\netstate(k)=j]{\m[][]{\Theta}(k+1)\m[T][]{\Theta}(k+1)}    \mathbreak 
= \expvalSym \left\{ \bar{\m{A}}(d(k)) \m[][]{\Theta(k)} \m[T][]{\Theta}(k) \m[T]{\bar{A}}(d(k)) \right. \mathbreak
\qquad \left. + \bar{\m{B}}\omega(k)\omega(k) \bar{\m[T]{B}} \middle| \netstate(k) = j \right\} \mathbreak 
 = \expvalSym \left\{ \bar{\m{A}}(d(k)) \m[][]{\Theta}(k) \m[T][]{\Theta}(k) \m[T]{\bar{A}}(d(k)) \middle| \netstate(k)=j \right\} \mathbreak
\qquad  + \bar{\m{B}}\sigma_{\omega}^{2} \bar{\m[T]{B}} . \label{eq:B2}
  \end{split}
\end{align}
Using the law of total expectation, we can write
\begin{align}
  \begin{split} &\expval[\netstate(k)=j]{\bar{\m{A}}(d(k))\m[][]{\Theta}(k)\m[T][]{\Theta}(k)\bar{\m[T]{A}}(d(k))} \mathbreak
     = \sum_{i \in \mathbb{B}} \expvalSym \Big\{\bar{\m{A}}(d(k))\m[][]{\Theta}(k)\m[T][]{\Theta}(k)\bar{\m[T]{A}}(d(k))  \mathbreak
 \qquad \Big| \netstate(k) = j, \netstate(k-1) = i \Big\} \mathbreak
     \times\prob[\netstate(k) = j]{\netstate(k-1)=i}.
  \end{split}
\end{align}
Using Bayes rule, this is rewritten to
\begin{align*} &\expval[\netstate(k)=j]{\bar{\m{A}}(d(k))\m[][]{\Theta}(k)\m[T][]{\Theta}(k)\m[T]{\bar{A}}(d(k))} \mathbreak 
                 = \sum_{i \in \mathbb{B}} \expvalSym \left\{\bar{\m{A}}(d(k))\m[][]{\Theta}(k)\m[T][]{\Theta}(k)\m[T]{\bar{A}}(d(k)) \right.\mathbreak
                 ~~~\left| \netstate(k) = j, \netstate(k-1) = i \right\} \mathbreak
                 ~~~\times \frac{\prob[\netstate(k-1) = i]{\netstate(k)=j}\prob{\netstate(k-1)=i}}{\prob{\netstate(k)=j}} \mathbreak
                 = \sum_{i \in \mathbb{B}} p_{ij} \expval[\netstate(k)=j,\netstate(k-1)=i]{\m[][]{\bar{A}}(d(k))} \mathbreak
                 ~~~ \times \expval[\netstate(k)= j,\netstate(k-1)=i]{\m[][]{\Theta}(k)\m[T][]{\Theta}(k)} \mathbreak
                 ~~~ \times \expval[\netstate(k)=j,\netstate(k-1)=i]{\m[T][]{\bar{A}}(d(k))} \mathbreak
                 ~~~ \times \frac{\prob{\netstate(k-1)= i}}{\prob{\netstate(k) = j}} \mathbreak
                = \sum_{i \in \mathbb{B}} p_{ij} \m[][j]{\mathcal{A}} \expval[\netstate(k-1)=i]{\m[][]{\Theta}(k)\m[T][]{\Theta}(k)}   \m[T][j]{\mathcal{A}}\mathbreak
~~~ \times\frac{\prob{\netstate(k-1)=        i}}{\prob{\netstate(k)= j}}.
\end{align*}
This can be inserted in (\ref{eq:B2}), such that
\begin{align}
  \begin{split}
    &\expval[\netstate(k)=j]{\m[][]{\Theta}(k+1)\m[T][]{\Theta}(k+1)} =    \mathbreak 
\sum_{i \in \mathbb{B}} p_{ij} \m[][j]{\mathcal{A}} \expval[\netstate(k-1)=i]{\m[][]{\Theta}(k) \m[T][]{\Theta}(k)} \m[T][j]{\mathcal{A}} \mathbreak 
\frac{\prob{\netstate(k-1)= i}}{\prob{\netstate(k)= j}} + \sigma_{\omega}^{2} \bar{\m{B}} \bar{\m[T]{B}}.
  \end{split}
\end{align}
Inserting this in (\ref{eq:BQ}) results in
\begin{align}
 &\m[][j]{Q}(k+1) =  \sum_{i \in \mathbb{B}} p_{ij} \m[][j]{\mathcal{A}} \expval[\netstate(k-1)=i]{\m[][]{\Theta}(k)\m[T][]{\Theta}(k)} \nonumber \mathbreak
\times \prob{\netstate(k-1)= i} \m[T][j]{\mathcal{A}} + \bar{\m{B}}\sigma_{\omega}^{2} \bar{\m[T]{B}}\prob{\netstate(k)=j} \nonumber \\
&= \sum_{i \in \mathbb{B}} p_{ij} \m[][j]{\mathcal{A}} \m[][i]{Q}(k) \m[T][j]{\mathcal{A}} + \sigma_{\omega}^{2}\bar{\m{B}} \bar{\m[T]{B}}\prob{\netstate(k)=j} \label{eq:Qk+1}
\end{align}

Since the \ac{NCS} \cref{eq:sysstability} is assumed to be \ac{MSS}, it is according to \cite[Theorem 3.33]{DT-MJLS} \ac{AWSS}. By defining $\m[][j]{Q} = \lim_{k \to \infty} \m[][j]{Q}(k), ~ j \in \mathbb{B}$, where every $\netstate(k)$ is aperiodic,
(\ref{eq:Qk+1}) becomes (\ref{eq:6}) and (\ref{eq:covQ}) whereas (\ref{eq:covQ}) becomes (\ref{eq:cov2state}).


\ifCLASSOPTIONcaptionsoff
  \newpage
\fi

\bibliographystyle{IEEEtran}

\bibliography{bare_jrnl}

\end{document}